\documentclass[12pt]{amsart}
\usepackage{amsmath,amsthm,amsfonts,amssymb}

\newtheorem{theorem}{Theorem}

\newtheorem{remark}{Remark}

\numberwithin{equation}{section}

\begin{document}

\title[Improved Uniqueness]
{An improved uniqueness result for a system of stochastic
differential equations related to the stochastic wave equation}
\author[C. Mueller, E. Neuman, M. Salins, and G. Truong]
{Carl Mueller
\and Eyal Neuman
\and Michael Salins
\and Giang Truong}

\address{Carl Mueller: Dept. of Mathematics
\\University of Rochester
\\Rochester, NY  14627}
\email{carl.e.mueller@rochester.edu}

\address{Eyal Neuman: Dept. of Mathematics
\\Imperial College London
\\ London, UK  SW7 2AZ}
\urladdr{http://eyaln13.wixsite.com/eyal-neuman}

\address{Michael Salins: Dept. of Mathematics and Statistics
\\Boston University
\\ Boston, MA 02215}
\urladdr{http://math.bu.edu/people/msalins/}

\address{Giang Truong:  Dept. of Mathematics
\\University of Rochester
\\Rochester, NY  14627}
\email{gtruong@u.rochester.edu}

\thanks{CM: Supported by a Simons grant.}

\begin{abstract}
We improve on the strong uniqueness results of \cite{glmns17}, which
deal with the following system of SDE.
\begin{align*}
dX_t&=Y_tdt  \\
dY_t&=|X_{t}|^{\alpha}dB_t
\end{align*}
and $X_0=x_0,Y_0=y_0$.  For $(x_0,y_0)\ne(0,0)$, we show that
short-time uniqueness holds for $\alpha>-1/2$.
\end{abstract}
\maketitle

\section{Introduction}
\label{section:introduction}

The purpose of this note is to improve a uniqueness result of
\cite{glmns17}.  First we state our result, and then we recall some
motivation.  Let $X_t,Y_t$ solve the following system of stochastic
differential equations (SDE) for $\alpha \in \mathbf{R}$.
\begin{align} \label{eq:sde}
dX_t&=Y_tdt  \\
dY_t&=|X_t|^{\alpha}dB_t,  \nonumber
\end{align}
with initial data $X_0=x_0,Y_0=y_0$.  Here $B_t$ is a standard
one-dimensional Brownian motion.
For the standard theory of SDE such as \eqref{eq:sde}, see Chapter V of \cite{protter05}.

We recall the results of Theorems 1.1 and 1.2 from \cite{glmns17}, which are stated together as follows.
\begin{theorem}[Gomez, Lee, Mueller, Neuman, and Salins] \label{th:ajcem}
If $\alpha>1/2$ and $(x_0,y_0)\ne(0,0)$, then (\ref{eq:sde}) has a
unique solution in the strong sense, up to the time $\tau$ at which the solution $(X_t,Y_t)$ first takes the value $(0,0)$ or blows up. Moreover the unique strong solution never reaches the origin.
\end{theorem}

In our main result we prove that the lower bound on $\alpha$ could be extended to $\alpha >-1/2$.
\begin{theorem} \label{th:main}
If $\alpha>-1/2$ and $(x_0,y_0)\ne(0,0)$, then (\ref{eq:sde}) has a
unique solution in the strong sense, up to the time $\tau$ at which the solution $(X_t,Y_t)$ first takes the value $(0,0)$ or blows up. Moreover the unique strong solution never reaches the origin.
\end{theorem}

 \begin{remark}
 The point $(x_0,y_0)=(0,0)$ plays a special role.  As proved
in Theorem 3 of \cite{glmns17}, if $0<\alpha<1$ then with this
initial condition both strong and weak uniqueness fail.
\end{remark}

Now we give some motivation for \eqref{eq:sde}.  Uniqueness questions
for SDE such as $dX_t=a(X_t)dt+b(X_t)dB_t$ have been studied for a long time.
Existence and uniqueness hold for Lipschitz coefficients $a,b$, see Section V.3 of
\cite{protter05}. The coefficient $a(x)$ can be badly behaved, but the
best result for $b(x)$, due to Yamada and Watanabe \cite{yw71}, is
that $b(x)$ should be H\"older continuous of order at least $1/2$.
However, Yamada and
Watanabe's method is essentially one dimensional, and does not carry
over to multidimensional systems except in special cases such  as radial symmetry.

For stochastic PDE, existence and uniqueness hold for most equations
in the case of Lipschitz continuous coefficients.  A case of special
interest is the SPDE for the superprocess,
\begin{equation*}
\partial_tu(t,x)=\Delta u(t,x)+|u(t,x)|^{1/2}\dot{W}(t,x), \quad x  \in\mathbf{R}, \, t\geq 0.
\end{equation*}
with appropriate initial data, usually nonnegative.  Here $\dot{W}(t,x)$ is two-parameter white noise.  For such initial data, weak uniqueness
among nonnegative solutions is known \cite{per02}, and strong uniqueness
among nonnegative solutions is an important unsolved problem.  If the
exponent $1/2$ is replaced by $\gamma>0$, then we know that strong
uniqueness holds among solutions taking values in $\mathbf{R}$ if
$\gamma>3/4$ \cite{mp11}, and both strong and weak uniqueness fail for $\gamma<3/4$
\cite{mmp14}.  The strong uniqueness results for $\gamma>3/4$ also
hold if $|u|^{\gamma}$ is replaced by a function of $u$ which is
H\"older continuous with index $\gamma$.

Much less is known about the stochastic wave equation
\begin{equation}
\label{eq:wave}
\partial_t^2u=\Delta u+|u|^{\alpha}\dot{W}(t,x)
\end{equation}
and analogous existence and uniqueness results are currently out of
reach.  Thus we are led to study SDE analogues of \eqref{eq:wave}
such as
\begin{equation*}
\ddot{u}(t)=|u(t)|^{\alpha}\dot{B}_t.
\end{equation*}
If we write $X_t=u(t)$ and $Y_t=\dot{u}(t)$, we arrive at \eqref{eq:sde}.

\section{Proof of Theorem \ref{th:main}}

First, recall that from Yamada and Watanabe \cite{yw71}, we know that the existence
of a weak solution together with strong uniqueness implies existence and uniqueness
in the strong sense.

\paragraph{\textbf{Step 1: Construction of a weak solution}}
When $0<\alpha \leq 1/2$, the construction of a weak solution and the proof that it almost surely never hits the origin is similar to the proof of Theorem 1.2 in \cite{glmns17}, hence it is omitted.

Assume now that $-1/2<\alpha \leq 0$, and fix the initial point $(x_0,y_0)\ne(0,0)$.

We use the following transformation which was used in the proof of Theorem 1.2 in \cite{glmns17}.
Define
\begin{equation}
\label{dife-h}
h(x):=\frac{1}{2\alpha+1}|x|^{2\alpha+1}\text{sgn$(x)$}, \quad h^{-1}(x):=(2\alpha+1)^{\frac{1}{2\alpha+1}}|x|^{\frac{1}{2\alpha+1}}\text{sgn$(x)$}.
\end{equation}
Observe that
$$
dh(x) = |x|^{2\alpha}dx, \quad dh^{-1}(x) = (2\alpha+1)^{\frac{-2\alpha}{2\alpha+1}}|x|^{-\frac{2\alpha}{2\alpha+1}}.
$$
Note that $h(x)$ is continuous and increasing in $\mathbf{R}$ even for
$-1/2<\alpha \leq 0$, and therefore the inverse function $h^{-1}(x)$ is well
defined. However, for $-1/2<\alpha < 0$, $dh(x)$ is infinite at the origin so the transformation in
Theorem 1.2 of \cite{glmns17} does not apply directly (see (3.4)--(3.6)
therein). Since $-1/2<\alpha \leq 0$, it follows that $dh^{-1}(x)$ is continuous in
$\mathbf{R}$.

Let
\begin{equation} \label{Y-V-eq}
\tilde{V}_t=h(x_0) + y_0 t + \int_{0}^{t}\tilde{B}_sds, \ \ \  \tilde{Y}_t = y_0 + \tilde{B}_t,
\end{equation}
where $\{\tilde{B}_t\}_{t\geq 0}$ is a standard Brownian motion.
We define the following time change
\begin{equation} \label{t-change}
 T(t) = \int_0^t(2\alpha + 1)^{-\frac{2\alpha}{2\alpha+1}} |\tilde{V}_s|^{-\frac{2\alpha}{2\alpha+1}} ds.
 \end{equation}
Note that
\begin{equation} \label{t-fin}
P(T(t) <\infty, \, \textrm{for all } 0 \leq t < \infty) =1,
\end{equation}
since $-\frac{2\alpha}{2\alpha+1} \geq 0$ for $-1/2<\alpha \leq 0$, and $\tilde{V}_s$ has continuous trajectories.

We further define the inverse time change,
\begin{equation}\label{inverse-t-change1}
T^{-1}(t) = \inf\{s\geq 0: T(s)>t\}.
\end{equation}
From Remark 5.2 in \cite{glmns17} we get that $|\tilde{V}_t| \vee |\tilde{Y}_t| \rightarrow \infty$ as $t\rightarrow \infty $, while both $\tilde{V}_t$ and $\tilde{Y}_t$ are recurrent process, hence it follows that  $\lim_{t\rightarrow \infty }T(t) = \infty$ a.s. and therefore
\begin{equation} \label{t-1-fin}
P(T^{-1}(t) < \infty, \, \textrm{ for all } 0\leq t< \infty) =1.
\end{equation}
Define
\begin{equation}  \label{x-t-def}
X_t = h^{-1}\big(\tilde{V}_{T^{-1}(t)}\big), \quad t\geq 0.
\end{equation}
First, we explicitly compute $T^{-1}(t)$:
\begin{align*}
  \frac{d}{dt} T^{-1}(t) &= \frac{1}{\frac{d}{ds} T(s)|_{s =T^{-1}(t)}} =(2\alpha + 1)^{\frac{2\alpha}{2\alpha+1}} |\tilde{V}_{T^{-1}(t)}|^{\frac{2\alpha}{2\alpha+1}} \\
  &  = (2\alpha + 1)^{\frac{2\alpha}{2\alpha+1}} |h(X_t)|^{\frac{2\alpha}{2\alpha+1}}
  = |h^{-1}(h(X_t))|^{2\alpha} = |X_t|^{2\alpha}.
\end{align*}
It follows that
\begin{equation} \label{t-1}
T^{-1}(t) = \int_0^t |X_s|^{2\alpha}ds.
\end{equation}
From (\ref{Y-V-eq}) and (\ref{t-1}) we get that
\begin{align*}
d\tilde{V}_{T^{-1}(t)} &=( y_0+\tilde B_{T^{-1}(t)}) dT^{-1}(t)  \\
&=( y_0+\tilde B_{T^{-1}(t)}) |X_t|^{2\alpha} dt.
\end{align*}
On the other hand, from (\ref{x-t-def}) we get,
\begin{align*}
d\tilde{V}_{T^{-1}(t)} &=dh(X_t) \\
&=  |X_t|^{2\alpha} dX_t.
\end{align*}
From (\ref{t-1-fin}) and (\ref{t-1}) it follows that the set $\{t\geq 0\, : \,   X_{t}=0\}$ has zero Lebesgue measure $P$-a.s. and therefore we have
\begin{equation}  \label{x-t}
dX_t = ( y_0+\tilde B_{T^{-1}(t)})dt.
\end{equation}
From (\ref{t-fin}) we have $\lim_{t\rightarrow \infty}T^{-1}(t) = \infty$, a.s., hence using (\ref{t-1}) we can define
\begin{equation} \label{y-t}
Y_{t} = y_0   + \tilde B_{T^{-1}(t)}.
\end{equation}
From the Dambis-Dubins-Schwarz theorem (see Revuz and Yor \cite{ry99},
page 181, Theorem 1.6) we get that $\{Y_{t}\}_{t\geq 0}$ satisfies
\begin{equation} \label{y-trans}
 Y_t = y_0   + \int_0^t|X_s|^{\alpha} d B^{(1)}_s,
\end{equation}
where $B^{(1)}_t$ is another standard Brownian motion.

From (\ref{Y-V-eq}) and (\ref{x-t})--(\ref{y-trans}) it follows that
\begin{equation}   \label{x-y-trans}
(X_t,Y_t) = (h^{-1}(\tilde{V}_{T^{-1}(t)}), \tilde{Y}_{T^{-1}(t)}),
\end{equation}
is a weak solution to (\ref{eq:sde}).

In was proved in Section 3 of \cite{glmns17} that $(\tilde{V}_t,\tilde{Y}_t)$ never equals $(0,0)$, that is,
\[P\big((\tilde{V}_t,\tilde{Y}_t)\ne(0,0)\text{ for } t>0\big)=1.\]
Together with (\ref{x-y-trans}) and (\ref{t-1-fin}) it follows that
\[P\big((X_t,Y_t)\ne(0,0)\text{ for } t>0\big)=1.\]

\paragraph{\textbf{Step 2: Proof of strong uniqueness}}

Let $(X^i_t,Y^i_t):i=1,2$ be two solutions of (\ref{eq:sde}) starting from $(x_0,y_0)\not = 0$, moreover let $\tau_n$ for a natural number $n$ be the first time $t$ at which either
\begin{equation*}
|(X^1_t,Y^1_t)|_{\ell^\infty}\wedge|(X^2_t,Y^2_t)|_{\ell^\infty}\leq2^{-n}\end{equation*}
or
\begin{equation*}
|(X^1_t,Y^1_t)|_{\ell^\infty}\vee|(X^2_t,Y^2_t)|_{\ell^\infty}\geq2^{n},
\end{equation*}
where $|(x,y)|_{\ell^\infty}=|x|\vee|y|$ is the $\ell^\infty$ norm.

Finally, as in the proof of Theorem 1.1 in \cite{glmns17},  let $Y^{i,n}_t = Y^i_{t \wedge \tau_n}$ and $X^{i,n}_t = \int_0^t Y^{i,n}_sds$.
Notice that $(X^{i,n},Y^{i,n})$ solve
\begin{align}
	\label{sys-trunc}
dX^{i,n}_t &= Y^{i,n}_tdt  \\
dY^{i,n}_t &= |X^{i,n}_t|^\alpha \mathbf{1}_{[0,\tau_n]}(t)dB_t   \nonumber
\end{align}
and that $(X^{i,n}_t, Y^{i,n}_t) = (X^i_t, Y^i_t)$ if $t \leq \tau_n$.
Define
\begin{equation*}
D_t=E\left[\left(X^{1,n}_t-X^{2,n}_t\right)^2\right].
\end{equation*}
Recall that $x \mapsto|x|^\alpha$ is a Lipschitz continuous function except in a neighborhood of $x=0$.
 As discussed in Section 2 of \cite{glmns17}, there is a sequence of stopping times
\begin{align}
  &\sigma^i_0 = 0 \nonumber\\
  &\sigma^i_{k+1} = \inf\{t >\sigma^i_k: X^{i,n}_t = 0\} . \nonumber
\end{align}
These stopping times form a discrete set and do not accumulate.

In order to prove uniqueness up to time $\tau_n$, it is enough to prove that $X^{1,n}_t = X^{2,n}_t$ for all $t \in [0,\sigma^i_k\wedge \tau_n]$ for any $k$ and any $i$. We do this in two steps.

First, assume that $x_0 \not =0$. We will argue that $\sigma^1_1 = \sigma^2_1$ and $X^{1,n}_t = X^{2,n}_t$ for all $t \in [0,\sigma^1_1]$. If $|X^{1,n}_t|\wedge |X^{2,n}_t| >0$ for all $t \in [0,\tau_n]$, then a minimum is attained and because the coefficients in \eqref{sys-trunc} are Lipschitz continuous when $|X^{i,n}_t|$ is bounded away from zero, standard uniqueness arguments show that $X^{1,n}_t = X^{2,n}_t$ for $t \in [0,\tau_n].$ So we assume that there exists $i \in \{1,2\}$ such that $\sigma^i_1\leq  \tau_n$. That is, at least one of the $X^{i,n}_t$ hits zero before $\tau_n$.
For $\delta<|x_0|$, let $\rho^\delta = \inf\{t>0: |X^{1,n}_t|\wedge |X^{2,n}_t| < \delta\}$. Because the coefficients of \eqref{sys-trunc} are Lipschitz continuous when $\delta< |X^{i,n}_t| $, standard arguments can be used to show that $X^{1,n}_t = X^{2,n}_t$ for all $t \in [0,\rho_\delta]$. By letting $\delta \to 0$ it is clear that $X^{1,n}_t = X^{2,n}_t$ for all $ t \in [0, \lim_{\delta \to 0} \rho_\delta)$. From the continuity of $X^{i,n}$, $i=1,2$ it follows that $\lim_{\delta \to 0} \rho_\delta = \sigma^1_1 \wedge \sigma^2_1$, the first time that one of the $X^{i,n}_t$ hits zero. Therefore, $X^{1,n}_t = X^{2,n}_t$ for all $t \in [0,\sigma^1_1\wedge \sigma^2_1)$ and by again by continuity we can conclude that $X^{1,n}_{\sigma^1_1\wedge \sigma^2_1} = X^{2,n}_{\sigma^1_1\wedge \sigma^2_1} =0$ so that $\sigma^1_1 = \sigma^2_1$.

Second, assume that $x_0=0$.

It is enough to prove
the uniqueness of the solutions to (\ref{sys-trunc}) starting at $X^{i,n}_0=0$ up to the first time that either one of $|X^{i,n}_t|$'s hits level $2^{-n}$. Therefore, we can restrict time $t$ to the interval $[0,\eta]$, where $\eta$ is
the first time $t<\tau_n$ at which
\begin{equation*}
|X^{1,n}_t|\vee |X^{2,n}_t|=2^{-n}.
\end{equation*}
If there is no such time, then let $\eta=\tau_n$.
Then using the strong Markov property we can restart the process at $\eta$ and use the previous step to prove uniqueness up to time $\sigma^1_1$.

Without loss of generality we can assume that  $y_0>0$.
Following the argument starting at the bottom of page 5 of \cite{glmns17}, we first
note that
\begin{equation*}
	X_t^{i,n}=\int_{0}^{t}\int_{0}^{s}|X^{i,n}_r|^\alpha \mathbf{1}_{[0,\tau_n]}(r)dB_rds.
\end{equation*}
By the Cauchy-Schwarz inequality and Ito's isometry, we get
\begin{align*}
	E\left[\left(X^{1,n}_t-X^{2,n}_t\right)^2\right]
	&\leq tE\int_{0}^{t}\left(\int_{0}^{s}
	\big(|X^{1,n}_r|^\alpha-|X^{2,n}_r|^\alpha\big)
	\mathbf{1}_{[0,\tau_n]}(r)dB_r\right)^2ds  \\
	&= tE\int_{0}^{t}\int_{0}^{s}
	\big(|X^{1,n}_r|^\alpha-|X^{2,n}_r|^\alpha\big)^2
	\mathbf{1}_{[0,\tau_n]}(r)drds    \\
	&\leq tE\int_{0}^{t}\int_{0}^{t}
	\big(|X^{1,n}_r|^\alpha-|X^{2,n}_r|^\alpha\big)^2drds  \\
	&\leq t^2E\int_{0}^{t}\big(|X^{1,n}_r|^\alpha-|X^{2,n}_r|^\alpha\big)^2dr.
\end{align*}
Thus, for the stopping time $\eta>0$ and any $t \in (0,\eta)$,
 \begin{align*}
D_t  \leq t^2E\int_{0}^{t}\big(|X^{1,n}_r|^\alpha-|X^{2,n}_r|^\alpha\big)^2dr.
\end{align*}

Now the mean value theorem gives, for $0<a<b$, that for some $c\in(a,b)$
we have
\begin{equation*}
b^\alpha-a^\alpha = \alpha c^{\alpha-1}(b-a) \leq |\alpha| a^{\alpha-1}(b-a).
\end{equation*}
Thus for $t\in[0,\eta]$, using the lower bound on $X^{i,n}_t$ from (2.3) in \cite{glmns17} we get
 \begin{equation}
\label{eq:new-gronwall}
D_t\leq |\alpha| 2^{-n(\alpha-1)} t^2\int_{0}^{t}r^{2\alpha-2}D_rdr.
\end{equation}
By assumption, for $i=1,2$,  $Y_t^{i,n}$
is almost surely continuous.
It follows that
\begin{equation} \label{y_0lim}
\lim_{t\downarrow0}\frac{X^{i,n}_t}{t}=
\lim_{t\downarrow0}\frac{1}{t}\int_{0}^{t}Y^{i,n}_rdr
=y_0
\end{equation}
exists.

From (\ref{sys-trunc}) we have
$$
|X^{i,n}_t| \leq \int_0^t|Y_s^{n,i}|ds \leq 2^nt,
$$
hence it follows that
$$
\left(X^{1,n}_t-X^{2,n}_t\right)^2 \leq 2^{2(n+1)}t^2.
$$
Then from dominated convergence we get
 \begin{equation*}
\lim_{t\downarrow0}\frac{D_t}{t^2} = (y_0-y_0)^2 = 0.
\end{equation*}

Let
\begin{equation*}
V_t=\frac{D_t}{t^2}.
\end{equation*}
By the above, $V_0=0$ exists as a limit.  Using \eqref{eq:new-gronwall}
we conclude
\begin{equation*}
V_t\leq C_n\int_{0}^{t}r^{2\alpha}V_rdr, \quad \textrm{for all } t \in (0,\eta),
\end{equation*}
and by Gronwall's lemma,
\begin{align*}
V_t&\leq V_0\exp\left(\int_{0}^{t}C_nr^{2\alpha}dr\right) \\
&\leq V_0\exp\left(\frac{C_n}{2\alpha+1}t^{2\alpha+1}\right)  \\
&=0.
\end{align*}

This shows uniqueness for $\alpha>-1/2$.

Finally, by using the strong Markov property and starting over at time $\sigma^1_1 = \sigma^2_1$, we can
extend our uniqueness result up to time $\sigma^1_2 = \sigma^2_2$.  By repeating this argument and using the fact that the $\sigma^i_k$ cannot accumulate, we can prove uniqueness up to time $\tau_n$.


\newcommand{\etalchar}[1]{$^{#1}$}
\providecommand{\bysame}{\leavevmode\hbox to3em{\hrulefill}\thinspace}
\providecommand{\MR}{\relax\ifhmode\unskip\space\fi MR }
\providecommand{\MRhref}[2]{%
  \href{http://www.ams.org/mathscinet-getitem?mr=#1}{#2}
}
\providecommand{\href}[2]{#2}

\end{document}